\documentclass[10pt]{IEEEtran}
\usepackage[dvips]{graphicx}
\usepackage{cite}
\usepackage{epsfig,enumitem}
\usepackage[colorlinks=true]{hyperref}
\usepackage{pstricks}
\usepackage{amsmath,amssymb,amsfonts}
\usepackage{color}
\usepackage{textcomp}
\def\BibTeX{{\rm B\kern-.05em{\sc i\kern-.025em b}\kern-.08em
    T\kern-.1667em\lower.7ex\hbox{E}\kern-.125emX}}

\newcommand{\dt}{\Delta t}

\begin{document}
\definecolor{lightgray}{gray}{0.9}

\title{Towards a comprehensive impossibility result for string stability}
\author{Arash Farnam \thanks{A.F. and A.S. are with the Department of Electronics and Information Systems (ELIS), Ghent University. Technologiepark 914, 9052 Zwijnaarde, Belgium.}, Alain Sarlette\thanks{A.S. is with the QUANTIC lab, INRIA Paris. 2 rue Simone Iff, 75012 Paris, France (email: alain.sarlette@inria.fr)}}
\date{\today}

\maketitle

\begin{abstract}
We provide a comprehensive impossibility result towards achieving string stability, i.e. keeping local relative errors in check with local controllers independently of the size of a chain of subsystems. We significantly extend existing results, from the LTI setting to any homogeneous controllers that can be nonlinear, time-varying, and locally communicating. We prove this impossibility for a set of definitions with various norm choices, including the $L_2$-type which is more standard in the literature and a BIBO type criterion. All results hold for a general discrete-time controller which should cover most applications.
\end{abstract}


\section{Introduction}\label{sec:intro}

String stability roughly requires that a chain of subsystems (e.g.~vehicles), controlled by local feedback from relative position measurements and subject to bounded local perturbations, keeps local deformations bounded independently of the size of the chain. In absence of further variables intervening in the system, local deformations are defined as deviations of the distances between consecutive vehicles from a constant target value. The topic has gathered attention from the observation that, for second-order integrators using \emph{any} LTI controller reacting to the distance to their predecessor, the amplification of some disturbance components along the chain is unavoidable \cite{1,2}. The historical application is a chain of acceleration-controlled vehicles, but other distributed systems like mechanical structures \cite{3,4}, or more fundamental models, should also benefit from string stability insight. This has initiated a rich line of results, among others: 
\begin{itemize}
\item Establishing various impossibility results and scalings for LTI controllers \cite{1,5,6}, mostly with subsystems reacting to their immediate predecessor \emph{and} immediate follower. Often PD controllers are used as a proxy for bandwidth limitation, but problems associated to integral control have also been characterized \cite{8}.
\item Designing a passive LTI controller, with symmetric coupling to predecessor and follower, to keep in check at least in some sense the effect of disturbances acting on the leading subsystem \cite{3}.
\item Adding a dependence on absolute velocity into the pure double integrator, to overcome the string instability issue \cite{10,12,13}. This can take the form of strong enough drag, which e.g.~for terrestrial vehicle applications are easy to ensure but would be in tradeoff with fuel efficiency; in other applications, ensuring strong drag may not always be trivial. Alternatively, the ``time-headway'' spacing policy makes the target inter-vehicle distance dependent on absolute velocity. For terrestrial vehicles again this appears easy to achieve, but questionable in terms of chain performance (possibly big spacing at high velocities, instead of moving efficiently like a big rigid body at all speeds); for other applications, the necessity to measure absolute velocity of each subsystem with respect to a common reference may even pose sensing questions.
\item Adding local communication between subsystems, as in adaptive cruise control \cite{14,15}. Remarkably, the solutions with this approach all incorporate as well a dependence on absolute velocity, as mentioned in the previous item. In the recent paper \cite{15a}, a variation on time headway is proposed for a generic controlled system and in the nominal case the controller uses both absolute state information and communication with the leader.
\end{itemize}
We can further refer the reader to the recent review \cite{MiddletonReview}. Research has also been carried out on properties related to, yet different from string stability. Without being comprehensive, we can mention the scaling of linear network eigenvalues with network size \cite{16}, the poor robustness of large linear networks under distributed sensing \cite{MyEELT}, and the coherence of linear networks in presence of stochastic noise \cite{18}.

This literature leaves several questions open. How much does each added element (dependence on absolute velocity, local communication capabilities, coupling symmetrically or asymmetrically to one or several vehicles ahead and/or behind) actually contribute towards solving the string stability issue? How much can be gained with more involved, e.g.~nonlinear controller designs? What is possible with digital, quantized controllers and realistic modern communication? How much is the LTI setting actually limiting?

The objective of the present note is to answer these questions towards a more precise understanding of string stability and related issues. As a small variation, we consider a discrete-time setting, which is closer to digital controller implementations and incorporates related ``natural'' constraints in a direct way; although probably most current controllers are digital and discrete-time anyways, the conclusions should carry over in practice to the continuous-time setting as we discuss below. We work only with \emph{relative} state measurements between neighboring subsystems, e.g.~we do not allow controller dependence on absolute velocity. In accordance with the existing literature, we consider an idealized model where input disturbances must be countered while assuming perfect measurements and communication. Our main result is:
\begin{itemize}
\item[] We establish that enabling nonlinear controllers, any couplings to a few vehicles in front and behind, any (e.g.~nonlinear, quantized, event-based) local communication, and controller dependence on the chain length, all together do not allow to design a controller which would achieve string stability with respect to bounded disturbances acting on the subsystems of the chain. We prove this for several variants of the string stability definition, and with as only main constraints: (i) the controller uses relative state measurements between neighboring vehicles; (ii) the controller is homogeneous, i.e. each vehicle (except boundaries) reacts in the same way to its neighbors; (iii) controller discretization step $\dt$ remains bounded away from zero despite increasing chain length.\newline
The proof comes down to identifying a general counterexample, then working out rather basic computations; we hence believe that the contribution mainly rests on the unprecedented generality of the conclusions. Essentially: \emph{string instability in a chain of second-order integrators is an unavoidable property of distributed sensing, for a (much) larger class of controllers than LTI.}
\end{itemize}

The note is organized as follows. Section \ref{sec:setting} describes the setting and several precise definitions of string stability. Section \ref{sec:impossibility} gives our main result, about a (very) general impossibility of string stability; it ends with a short discussion of discrete-time vs.~continuous-time controllers, and an illustrative simulation of the system behavior under a ``bad'' disturbance. Section \ref{sec:conclusion} concludes the note with an outlook on remaining open points and possible further implications.


\section{Problem definition}\label{sec:setting}


\subsection{Model setting}

We consider a chain of undamped second-order integrators, called ``vehicles'' for concreteness but which may represent any other (typically mechanical) subsystems coupled with local interactions. In the literature, this continuous-time system
$$\ddot{x}_k(t) = u_k(t)$$
with $x_k$ the position and $u_k$ the input signal of subsystem $k$ at time $t$, is the standard starting point for string (in)stability \cite{1,2,3,4,5,6,MiddletonReview}, with continuous-time e.g.~LTI controllers defining the $u_k$. We here consider its (exact) equivalent under discrete-time control. We assume that each subsystem is controlled with a discrete-time control logic, computing at each time $t = n\, \dt$ with $n \in \mathbb{Z}$, a control signal that will be applied during the whole interval $(t,t+\dt]$ as input to each double integrator; thus $\dt$ is the time increment of the discrete-time controller. We can then integrate exactly the second-order dynamics over the time interval $(t,t+\dt]$ to obtain:
\begin{eqnarray}
\label{eq:model} v_k(t+\dt) & = & v_k(t) + u_{k,1}(t) + d_{k,1}(t)\\
\nonumber x_k(t+\dt) & = & x_k(t) + v_k(t)\, \dt + u_{k,2}(t) + d_{k,2}(t) \; ,
\end{eqnarray}
where $x_k,v_k \in \mathbb{R}$ for $k=0,1,2,...,N$ denote position and velocity respectively; the control inputs $u_{k,1},u_{k,2}$ result from integrating the control signal applied during the interval $(t,t+\dt]$ respectively once and twice; and the disturbances $d_{k,1},d_{k,2}$ result from similarly integrating the continuous-time input disturbances. The key point in calling \eqref{eq:model} an \emph{exact} discrete-time equivalent is that $u_k(t)$ over the interval $(t,t+\dt]$ does not depend on the state values after time $t$.
 
A key constraint is that the feedback signals $u_{k,1},u_{k,2}$ applied between $t$ and $t+\dt$ must be based on only \emph{relative} displacement measurements $e_k = x_{k-1}-x_k$ between consecutive vehicles in the chain, including possibly relative speeds $\dot{e}_k = v_{k-1} - v_k$, taken at time $t$. We then consider the following general controller, where $e_{k:\ell}$ denotes for $\ell > k$ the set of values $e_k,e_{k+1},...,e_\ell$:
\begin{eqnarray}
\label{eq:controller} u_{k,1} & = & f_1(e_{(k-m_1):(k+m_2)},\dot{e}_{(k-m_1):(k+m_2)},\\
\nonumber& & \phantom{KKKKKKKK} c_{k,+},c_{k,-},\xi_k,N,t) \; ,\\
\nonumber u_{k,2} & = & f_2(e_{(k-m_1):(k+m_2)},\dot{e}_{(k-m_1):(k+m_2)},\\
\nonumber& & \phantom{KKKKKKKK} c_{k,+},c_{k,-},\xi_k,N,t) \; ,\\
\nonumber c_{k+1,+} & = & g_1(e_{(k-m_1):(k+m_2)},\dot{e}_{(k-m_1):(k+m_2)},\\
\nonumber& & \phantom{KKKKKKKK} c_{k,+},c_{k,-},\xi_k,N,t) \; , \\
\nonumber c_{k-1,-} & = & g_2(e_{(k-m_1):(k+m_2)},\dot{e}_{(k-m_1):(k+m_2)},\\
\nonumber& & \phantom{KKKKKKKK} c_{k,+},c_{k,-},\xi_k,N,t) \; , \\
\nonumber \xi_k(t+\dt) & = & h(\xi_k(t),e_{(k-m_1):(k+m_2)},\\
\nonumber&& \phantom{aaaaaa} \dot{e}_{(k-m_1):(k+m_2)},c_{k,+},c_{k,-},N,t) \, .
\end{eqnarray}
Here $m_1$ (resp.~$m_2$) is a finite number of agents ahead (resp.~behind); $c_{k,+},c_{k,-} \in \mathbb{R}^{n_c}$ with $n_c$ some bounded integer are communication signals from $\{k-m_1,...,k-1\}$ to $k$ and from $\{k+1,...,k+m_2\}$ to $k$ respectively; the $\xi_k \in \mathbb{R}^{n_{\xi}}$ for some finite integer $n_\xi$ allow for dynamical controllers with finite memory; and $f_1,f_2,g_1,g_2,h$ are arbitrary functions, with minimal regularity just to ensure that the solution to the dynamical system is well-defined at all times. In particular, by specifying a particular profile of $u(\tau)$ over the discretization time span $\tau \in [t,t+\dt)$ between two measurement updates, one can command $f_1=\int_{t}^{t+\dt} u(\tau) d\tau$ and $f_2 = \int\int_{t}^{t+\dt} u(\tau) d\tau^2\,$ independently. The controllers \eqref{eq:controller} are applied by all vehicles $k \in (m_1,N-m_2)$, whereas adapted versions are applied by the $m_1$ leading and $m_2$ last vehicles. The adapted versions will play no role in the proof. 

Note that $m_1$ and $m_2$ are just opportunities, without obligations, for a vehicle to take into account what happens several vehicles ahead and/or behind. The setting includes \emph{unidirectional} chains $m_2=0$, where vehicles only get inputs from predecessors; in particular, $m_2=0$ and $m_1=1$ corresponds to reacting just to the preceding vehicle. It also includes the bidirectional \emph{symmetric} chain, like encountered in homogeneous mechanical coupling \cite{3,4}, by imposing front-to-back symmetry on \eqref{eq:controller}.

Our aim here is to leave as much freedom as possible on controller design. In particular, the functions $f_1,f_2,g_1,g_2,h$ can be nonlinear and time-dependent (e.g.~modulated at specific frequencies), thereby vastly extending the traditional LTI setting. Our only true restrictions on control design are: 
\vspace{2mm}

\noindent \textbf{Assumptions implied by the model:}
\begin{itemize}
\item \emph{The controllers must be homogeneous along the chain,} i.e.~the functions $f_1,f_2,g_1,g_2,h$ do not depend on vehicle index $k$ and the internal variables are initialized with the same default values for each $k$. 
\item \emph{The digital controller has a finite update time $\dt$}, fixed independently of chain length $N$.
\item \emph{The control commands are based on \emph{relative} state information only}, i.e.~there is no dependence on absolute states that would be obtained with respect to some common reference, like absolute velocity.
\end{itemize}
The first point can be relaxed to approximately homogeneous in our proof; future work, if deemed relevant, could address the fully heterogeneous setting e.g.~with adversarial noise reshaping on the basis of Sec.\ref{sec:impossibility}. The second point is discussed after the main result. Regarding the third point, we recall that the case with control using \emph{absolute} velocity has been solved \cite{10,12}. It is not clear however whether relying solely on absolute velocity for achieving string stability is acceptable in all practical contexts. This is also the original academic question: string stability without using any absolute information.


\subsection{Control objective}

The goal of the controller is to achieve \emph{string stability}. Roughly said, this means avoiding that the performance gets unboundedly worse when the chain length $N$ grows. Formally, several definitions have been proposed \cite{1,2,3,6,8,15,MiddletonReview}. Their common point is to focus on stabilizing the $e_k$, i.e.~the distances between consecutive vehicles. This makes sense e.g.~for collision-avoidance, for maintaining a tight platoon, or to avoid too extreme accelerations on the last vehicles of the chain, and is a weaker requirement than controlling the position of every vehicle with respect to the leader. Definitions vary in the way they consider the distribution over vehicles, as a sum or individually, the disturbances, as signals or initial conditions, and the norm over time of the signals \cite{15}; see e.g.~\cite{MiddletonReview} for a review. We here translate the definitions in presence of input disturbances to discrete-time control.\vspace{2mm}

\noindent \textbf{Definition 1:} \emph{For positive integers $p,q$, the $\ell_{p,q}$ string stability requires that there exist $C_1,C_2>0$ independent of $N$ such that: for any disturbances satisfying} 
\vspace{2mm}

$\,
\sum_{s=1,2}\; \sum_{k=0}^N \left(\sum_{n \in \mathbb{Z}} |d_{k,s}(n\, \dt)\tfrac{1}{\dt^s}|^p \dt\,\right)^{q/p}  < C_1 \; ,
$
\vspace{2mm}

\noindent \emph{it is ensured that $\; \sum_{n \in \mathbb{Z}} |e_j(n\, \dt)|^p \dt < C_2$ for each $j$.}\vspace{2mm}

\noindent \textbf{Definition 2:} \emph{For positive integers $p,q$, the $(\ell_p,\ell_q)$ string stability requires that there exist $C_1,C_2>0$ independent of $N$ such that: for any disturbances satisfying}
\vspace{2mm}

$\,\sum_{s=1,2}\; \sum_{k=0}^N \left(\sum_{n \in \mathbb{Z}} |d_{k,s}(n\, \dt)\tfrac{1}{\dt^s}|^p \dt\,\right)^{q/p} < C_1 \; ,
$
\vspace{2mm}

\noindent \emph{it is ensured that $\; \sum_{k=0}^N \left(\sum_{n \in \mathbb{Z}} |e_k(n\, \dt)|^p \dt \right)^{q/p}  < C_2$.}\vspace{2mm}

\noindent \textbf{Definition 3:} \emph{For positive integer $p$, the $(\ell_p,\ell_\infty)$ string stability requires that there exist $C_1,C_2>0$  independent of $N$ such that: for any disturbances satisfying}
\vspace{2mm}
 
$\sum_{n \in \mathbb{Z}} |d_{k,s}(n\,\dt)\tfrac{1}{\dt^s}|^p \dt <C_1\,,\;\; s=1,2\; \text{ for all } k\,,$
\vspace{2mm}

\noindent \emph{it is ensured that $\sum_{n \in \mathbb{Z}} |e_j(n\,\dt)|^p \dt < C_2$ for all $j$.}
\vspace{2mm}

\noindent \textbf{Definition 4:} \emph{The $(\ell_\infty,\ell_\infty)$ string stability requires that there exist $C_1,C_2>0$  independent of $N$ such that: for any disturbances satisfying $\; |d_{k,s}(t)|/\dt^s<C_1 \;$ for $s=1,2$ and all $k,t$, it is ensured that $\; |e_j(t)| < C_2\; $ for all $j,t$.} \hfill $\square$
\vspace{2mm}

The $\dt$ factors are introduced as $d_{k,1}$ and $d_{k,2}$ are supposed to result from integrating (respectively once and twice) a continuous-time signal $d_k(t)$ over the time interval $(t,t+\dt]$. In principle we could just choose units such that $\dt=1$, but we keep $\dt$ for later discussion. The main point about string stability is to satisfy the constraints with constants \emph{independent of the number of vehicles $N$.}

A priori, Definition 1 is the weakest since it imposes a bound on the \emph{vector} norm of disturbance inputs,  but in return it only requires each individual $e_j$ to have a bounded signal norm. While this may appear quite asymmetric, it is a necessary condition for achieving stronger versions of string stability, and it has been considered a lot in the literature; e.g.~until recently this was the only proven working definition when adding a time-headway policy (i.e.~a controller depending on \emph{absolute velocity}, as we exclude here \cite{10,12,13}). The other defnitions each consider the same norm on input disturbances and on output errors, either summing them over vehicles (Def.2) or not (Def.3, Def.4). The most popular norm has been $p=q=2$, especially the time-integration with a $2$-norm has attracted a lot of attention thanks to its equivalent formulation in frequency domain. Definition 4 changes the treatment of time to formulate a BIBO type version of string stability. It seems to have attracted little attention in the literature, probably because most of the literature has focused on LTI frequency domain approaches. However, we would argue that this version is closest to practical concerns; and, as we will show that the string instability issue is not limited to the LTI setting, we will use tools that can treat this definition explicitly too. In the particular context of \eqref{eq:model},\eqref{eq:controller}, or when moreover e.g.~assuming controllers to be LTI, it might well be that some of the above definitions hold strictly together; in absence of further evidence, we will treat them all.\vspace{2mm}

\noindent \textbf{Remark 1:} The definitions were initially stated in the linear context, where $C_1$ and $C_2$ can be rescaled such that it makes no difference in which order they are chosen (e.g.~variants like ``for each $C_1$, there exists a $C_2$'' become equivalent to our statement). In the nonlinear context this might differ, and we have chosen the weaker constraint: thus our impossibility results will also hold for stronger variants.\hfill $\square$\\ 

\noindent \textbf{Remark 2:} The definitions require to check system behavior under the \emph{worst disturbance} satisfying theconstraints. This may seem quite natural when checking stability or disturbance rejection. We will see though that it appears quite demanding for string stability. Some authors were able to obtain more positive results by considering disturbances restricted to the leading subsystem only, see e.g.~\cite{3}, or assuming a probability distribution on input disturbances, see e.g.~\cite{18}. We here stick to the most standard definition; after giving our results, we will come back to discuss this choice.\hfill $\square$


\section{Main Impossibility Result}\label{sec:impossibility}

We now prove that Definitions 1-4 are all impossible to satisfy even with a general controller as allowed by \eqref{eq:controller}. The main idea of the proof is to construct a disturbance input which is badly countered by any distributed controller. While exactly solvable situations are hard to find, we take advantage of a simple construction that focuses on the central part of the chain only, in order to give a lower bound on the induced error.

A simulation illustrating the behavior of the system under this construction can be found in Section \ref{sec:simulation}.


\subsection{A badly countered disturbance situation}\label{sec:badpert}

Consider disturbances of the following form:
\begin{eqnarray}
\nonumber d_{k,1}(t) = d_{k,2}(t) = 0 && \text{for all } t<0,\;\; k=0,1,...,N ;\\
\label{eq:disturbance} 
\left.\begin{array}{r}
d_{k,1}(t) = \tfrac{\alpha k \dt}{N},\\
d_{k,2}(t) = \tfrac{\alpha k \dt^2}{N} 
\end{array}\right\rbrace
&& \begin{minipage}{0.2\textwidth} 
for all  $t=0,\dt,...,T$, $k=0,1,...,N \; ;$
	\end{minipage}\\
\nonumber d_{k,1}(t) = d_{k,2}(t) = 0 && \text{for all } t>T,\;\; k=0,1,...,N \; ,
\end{eqnarray}
with constants $\alpha>0$ and $T>0$ to be specified later.

To compute the evolution of the system under these disturbances, the trick is to exploit the finite propagation speed of signals along the chain --- namely at most $(m_1,m_2)$ vehicles per time step --- in order to restrict our attention to a central subset of vehicles, for which the computations are easy.\vspace{2mm}

$\bullet$ Consider the evolution of $e_k$ and $\dot{e}_k$ over one time step, when the $N+1$ vehicles all start with the same state $x_k(0)=v_k(0)=0$ for all $k$, and with controllers initialized at $\xi_k=c_{k,+}=c_{k,-}=0$ for all $k$. We get
\begin{eqnarray*}
e_{k}(\dt) & = & e_k(0) + \dt \, \dot{e}_k(0) + u_{k-1,2}(0)-u_{k,2}(0)\\
&& \;\;\;  + d_{k-1,2}(0) - d_{k,2}(0)\\
& = &  u_{k-1,2}(0)-u_{k,2}(0) + \alpha \dt^2/N \; ; \\
\dot{e}_{k}(\dt) & = & \dot{e}_k(0) + u_{k-1,1}(0)-u_{k,1}(0)\\
&& \;\;\; + d_{k-1,1}(0) - d_{k,1}(0)\\
& = &  u_{k-1,1}(0)-u_{k,1}(0) + \alpha \dt/N \; .
\end{eqnarray*}
Since the $e_k(0)$ and $\dot{e}_k(0)$ are all equal, it is clear that the control inputs are all equal too, i.e.~$u_{k-1,1}=u_{k,1}$ and $u_{k-1,2}=u_{k,2}$, at least for all vehicles with $m_1<k<N-m_2$. For those vehicles, \emph{completely irrespectively of the controller chosen}, we have
\begin{align*}
e_k(\dt) = \alpha \dt^2/N \text{ and } \dot{e}_k(\dt) = \alpha \dt/N \;, \\ \text{ for all } m_1<k<N-m_2 \, .
\end{align*}
Also the $c_{k,+}(\dt),c_{k,-}(\dt)$ and $\xi_k(\dt)$ will be equal.\vspace{2mm}

$\bullet$ Now consider a time $t=n\dt$ for some integer $n>0$ and assume that all the state variables satisfy equalities $e_k(t)=e_j(t)$, $\dot{e}_k(t)=\dot{e}_j(t)$, $c_{k,+}(t)=c_{j,+}(t)$, $c_{k,-}(t)=c_{j,-}(t)$ and $\xi_k(t)=\xi_j(t)$ for all $j,k \in [N_\text{lead},\;N-N_\text{tail}]$, for some integers $N_\text{lead},N_\text{tail}>0$. Slightly extending the above example, we get:
\begin{eqnarray}
\label{eq:as1} e_{k}(t+\dt) & = &  e_k(t) + \dt \, \dot{e}_k(t) + u_{k-1,2}(t)-u_{k,2}(t)\\
\nonumber & & \;\; +\; d_{k-1,2}(t) - d_{k,2}(t)\\
\nonumber & \hspace{-15mm} = & \hspace{-8mm} e_k(t) + \dt \, \dot{e}_k(t) + \alpha \dt^2/N \;=\; e_j(t+\dt) \\
\nonumber & \hspace{-15mm} & \hspace{-13mm} \text{for all } j,k \in [N_\text{lead}+m_1,\;N-(N_\text{tail}+m_2)]\; ;\\
\nonumber \dot{e}_{k}(t+\dt)& = & \dot{e}_k(t) + \alpha \dt/N \;=\; \dot{e}_{j}(t+\dt) \\
\nonumber & \hspace{-15mm} & \hspace{-8mm} \text{ for all } j,k \in [N_\text{lead}+m_1,\;N-(N_\text{tail}+m_2)]\; ,
\end{eqnarray}
and similarly we maintain $c_{k,+}(t)=c_{j,+}(t)$, $c_{k,-}(t)=c_{j,-}(t)$ and $\xi_k(t)=\xi_j(t)$ for that subset of vehicles.\\

By iterating this argument we get the following property.\vspace{2mm}

\noindent \textbf{Lemma 1:} \emph{Consider the system \eqref{eq:model},\eqref{eq:controller} subject to the particular disturbance \eqref{eq:disturbance} and zero initial conditions. Then for any (well-defined) controller choice, the solution satisfies:
\begin{eqnarray}\label{eq:solution}
e_k(t) &=&  t(t+\dt) \; \alpha\;/\;(2N)\\
\nonumber \dot{e}_k(t) &=& t\; \alpha/N \, ,
\end{eqnarray}
for all $t \in [0,T]$ and all $k \in (\frac{t}{\dt}\,m_1,\;N-\frac{t}{\dt}\,m_2)$.}\vspace{2mm}

\noindent \underline{Proof:} The main argument is provided by the explanations preceding the statement. From \eqref{eq:as1}, the $\dot{e}_k$ is obtained as a sum of $t/\dt$ times the bias $\alpha \dt/N$. Then replacing this into the expression of $e_k$ in \eqref{eq:as1}, one observes that the increment of $e_k$ at time $n=t/\dt$ is linear in $n$, so the standard formula for a linearly progressing series gives the result. \hfill $\square$\\

To be useful at time $t$, the solution \eqref{eq:solution} of Lemma 1 should cover at least $1$ vehicle, i.e.~$N-\frac{t}{dt}(m_1+m_2)\geq 1$. For fixed $m_1, m_2$ and $dt$, we can ensure to have a valid solution for at least $N/2$ vehicles over the interval $[0,t]$, by taking $t=\frac{N \, dt}{2(m_1+m_2)}$. The shortest disturbance that will lead to the result of Lemma 1 on $N/2$ vehicles for $t=\frac{N \, dt}{2(m_1+m_2)}$, is by taking $T=\frac{N \, dt}{2(m_1+m_2)}$. As $m_1$ and $m_2$ are constants independent of $N$, we essentially suggest to select the duration $T$ of the ``bad'' input disturbance to be of order $N \, dt$.


\subsection{Consequences for string stability}

We now investigate what the above construction implies for string stability. First take Definition 1. For the proposed disturbance, the condition on $d_{k,1},\;d_{k,2}$ becomes:
$$2 T^{q/p}\, (\tfrac{\alpha}{N})^q\, {\textstyle \sum_{k=0}^N}\; k^q \; < C_1 \; .$$
For large $N$, the dominating term in the sum is $N^{q+1}$ such that we need in fact $2 T^{q/p} \alpha^q N < C_1$ or in other words, $\alpha$ of order $1\, /\, (N^{1/q}T^{1/p})$. This fixes the allowed disturbance amplitude as a function of its duration and of $N$.

For the vehicles covered by Lemma 1, we then have
\begin{eqnarray*}
\sum_{n\in\mathbb{Z}} |e_j(n\,\dt)|^p \dt & \geq & \sum_{n=0}^{T/\dt} |e_j(n\,\dt)|^p \,\dt\\
& \simeq & (\tfrac{\alpha}{2N})^p \sum_{n=0}^{T/\dt}  (n \dt)^{2p} \, \dt \\
& \simeq & (\tfrac{\alpha}{2N})^p \, \dt^{2p+1} \, (\tfrac{T}{\dt})^{2p+1} \; .
\end{eqnarray*}
Towards the last line we have again used the dominating term in the sum. Combining this with the just obtained bound on $\alpha$ and with taking $T$ of order $N\, \dt$ as suggested at the end of Section \ref{sec:badpert}, we obtain that $\sum_{n\in\mathbb{Z}}\, |e_j(n\,\dt)|^p$ is at least of order $\; N^{p-p/q}\, \dt^{2p}\;$.

A similar argument can be repeated for the other definitions of Section \ref{sec:setting}, yielding the following results.\vspace{2mm}

\noindent \textbf{Theorem 2:} \emph{For the system \eqref{eq:model},\eqref{eq:controller}, there exist disturbances $d_{k,1}$ and $d_{k,2}$ satisfying the required respective bounds according to the definitions of Section \ref{sec:setting} and such that, irrespectively of any (well-defined) controller choice, for large $N$:\vspace{2mm}
\newline \textbf{[Definition 1]:} $\sum_{n} |e_j(n\,\dt)|^p$ grows as $\,N^{p-p/q}\, \dt^{2p}\,$;\vspace{2mm}
\newline \textbf{[Definition 2]:} $\sum_{k=0}^N \left( \sum_{n} |e_k(n\, \dt)|^p \right)^{q/p}$ grows as $N^q \dt^{2q}$;\vspace{2mm}
\newline \textbf{[Definition 3]:} $\sum_{n} |e_j(n\, dt)|^p$ grows as $N^p \dt^{2p}$;\vspace{2mm}
\newline \textbf{[Definition 4]:} $|e_j(t)|$ grows as $N \dt^2$.}
\vspace{2mm}

\noindent \underline{Proof:} We will always assume $T$ of order $N \dt$ and consider the output errors for $t\in [0,T]$. The computation for Def.1 is given above. For Def.2 it is the same, but taking the power $q/p$ and summing the disturbance over the number of vehicles for which Lemma 1 is valid -- this can be of order $N$ as mentioned in the last paragraph of Section \ref{sec:badpert}. For Def.3, the disturbance can be larger i.e.~$\alpha$ of order $1/T^{1/p}$, and with respect to the Def.1 computation this adds a factor $N^{p/q}$ to the output error. For Def.4, we can have $\alpha$ of order 1, and since the result of Lemma 1 is valid for $T$ of order $N \dt$ we can have $e_k(t)$ of order  $T(T+1) \; \alpha/N \sim N^2 \dt^2/N$. \hfill $\square$\\

For $\dt$ fixed and $N$ growing to infinity, this result establishes impossibility to satisfy any of the definitions of string stability given in Section \ref{sec:setting}, except possibly Def.1 with $q=1$ (which does not appear to have any practical significance, see comments about the Definitions in Section \ref{sec:setting}). This impossibility is established in a very general setting, allowing unidirectional or bidirectional symmetric or asymmetric coupling, looking a number of vehicles ahead and behind (as long as that number is independent of $N$), communicating with neighbors with any encoding/decoding schemes with possibly packets and event-based logic, and processing all this in an arbitrary nonlinear control system with memory. In particular, even perfect local communication among the vehicles is not sufficient, on its own, to ensure string stability. It is thus no wonder that vehicle chain controllers with realistic communication channels have so far required an additional feedback from absolute velocity to achieve string stability \cite{14,15,15a}.


\subsection{How telling is the discrete-time controller setting?}

The reader will have noticed that the above impossibility breaks down if we let $\dt$ converge to zero fast enough as $N$ grows to infinity. While this does not look like a practical solution, it may express a relevant tradeoff; and, it does create a gap with the pure continuous-time literature. We will thus briefly comment on the comparison of this result with the literature on continuous-time, typically LTI systems.

As a common point with the literature, low-frequency disturbances indeed appear to cause most of the problem in continuous-time string instability proofs for LTI systems. For those cases, e.g.~PD coupling with nearest neighbors, adding the contributions of all the neglected vehicles and time-steps to the norm of $e$ would yield string instability in continuous-time too. In this sense, we can expect that the arguments leading to Thm.2 are too optimistic in the sense that simply letting $\dt \rightarrow 0$ would not actually solve the issue (see also the simulations below). More generally, as our result only follows a sufficient construction, Theorem 2 proves that it is \emph{necessary} -- yet possibly not even sufficient -- to let $\dt$ go to zero with increasing $N$ in order to satisfy string stability. So let us try to list and discuss which controller features would typically go with a very small $\dt$:
\begin{itemize}
\item One obvious effect of smaller $\dt$ is faster communication across the vehicle chain. If one could communicate arbitrarily fast, perfectly and without measurement errors, then each vehicle $k$ could get very fast knowledge of $e_1+e_2+...+e_k = x_k - x_0$. One can then obviously achieve string stability: just control each $x_k-x_0$ independently to stabilize each vehicle with respect to the leader. The ``distributed system'' setting and chain size $N$ play no role anymore. Of course this idealized situation is unrealistic. In reality, precision of a message (and in fact of a measurement) is in a clear tradeoff to update speed. As long as the communication bandwidth per signal remains bounded when $N$ increases, the imperfections resulting from smaller $\dt$ are likely to counterbalance the apparent benefits of smaller $\dt$ from our perfect-communication model.

\item Setting sensing and communication aside, in practice, the controller's discretization step $\dt$ is chosen as the desired dwell-time before vehicle $k$ reacts to a measurement; thus in practice $\dt$ converging to zero would mean, controller bandwidth tending to infinity, pointing towards controllers with gain increasing as a function of $N$. It is known indeed that academically speaking, this can provide string stability: in continuous-time, without communication, a LTI controller whose gain increases fast enough with $N$, can ensure string stability. However, as the control gain keeps increasing towards infinity, effects of unmodeled system limitations and imperfections cannot be neglected forever and practical problems are likely to appear.

\item Theorem 2 thus shows anyways that string stability is, at best, not robust to time-discretization. This is important to know towards system simulations, where situations that work only for infinitesimal $\dt$ are quickly considered non-robust for all practical purposes. In a sense, testing robustness to finite $\dt$ can even be mathematically compared to the traditional requirement of ``no poles cancellation'' in the continuous-time setting. Indeed, allowing a decreasingly small $\dt$ without any measurement noises can be compared to allowing precise computation of $\lim_{\dt \rightarrow 0} \frac{s(t+\dt)-s(t)}{\dt}$ for a signal $s$, i.e.~evaluating pure derivatives. For a double-integrator, this implies the possibility of pole cancellation at zero frequency, which is almost always excluded.

\item The dependence on $\dt$ is rooted in the fact that we analyze the system before the signals from the edges of the chain reach all the vehicles and make a detailed analysis harder. This does not mean of course that the vehicle chain would automatically be stabilized as soon as the signals from the edges have crossed the chain, see e.g.~the simulations in Section \ref{sec:simulation}. In this sense, it appears that the boundary controllers would play a key role towards string stabilizing the system with infinitesimal $\dt$, similarly to PDE control.
\end{itemize}

These arguments give strong indications to conjecture that string stability would be impossible with any ``reasonable'' homogeneous, possibly nonlinear, and communicating controllers, in continuous-time too. At this point of detail, we might argue as well that the digital-controller model is in fact closer to applications, than the traditional continuous-time one.
\vspace{3mm}

\noindent \textbf{Remark 3:} To further connect this result to existing work, we can look at how the chain reacts to disturbances acting on the first subsystem only. This has indeed been considered in several continuous-time LTI studies, which we first review now. For unidirectional coupling, reaction to leader-disturbances is sometimes viewed as a major indicator of general behavior \cite{1,2}. For bidirectional chains under symmetric coupling, although impossibility results are known from e.g.~\cite{8}, under the condition of disturbance restricted to the leading vehicle, $\ell_{2,2}$ string stability has been established in \cite{3}. Pushing further the idea of \cite{16} and of \cite{20,20b} about possible advantages of slight mistuning in the controller symmetry, we have proved in \cite{21} that a sufficiently asymmetric bidirectional PD controller is sufficient to ensure also the stronger versions of string stability, with respect to disturbances restricted to a fixed number of leading vehicles. The proof uses an analytic almost-inversion of the system equations based on forward and backward flows, loosely inspired from \cite{23}.

This line of work can be related in two ways to the result of the present paper. First, we have checked that a discrete-time version of our string stability result in \cite{23} can be worked out perfectly well for the model \eqref{eq:model},\eqref{eq:controller}, see \cite{PhDThesis}. This suggests that our discrete-time model does enable positive results when the continuous-time model does, i.e.~it adds evidence in favor of \eqref{eq:model},\eqref{eq:controller} not being essentially more constraining than the more standard continuous-time approach. Second, the disturbance proposed in Section \ref{sec:badpert} is extensively distributed along the vehicle chain, as a function of $N$. This is consistent with a picture of two regimes: when disturbances act on a few vehicles (at known places!), it may be possible to reject them in a string stable way; however when they are distributed along the whole chain, there is no way to achieve string stability on the basis of relative measurements only. \hfill $\square$

\subsection{Illustrative simulation}\label{sec:simulation}.

We can of course only illustrate the string instability for a particular choice of controller. However, trusting in the simple analysis of Section \ref{sec:badpert}, our main argument is independent of controller choice. Thus, we will just show how indeed the solution of Lemma 1 appears for a simple linear controller without communication. We choose this simplicity to avoid selecting too many elements in the controller design ``arbitrarily'' --- since according to Theorem 2, any attempt is anyways doomed to fail. Complementarily the simulation shows what happens once the ``boundary effects'' have propagated throughout the chain i.e.~when the solution of Lemma 1 is not valid anymore; this \emph{will} depend on the choice of controller, but it falls outside the scope of the present paper. We thus suggest that the reader should not draw too strong conclusions from what happens outside the scope of Lemma 1 with this particular controller.

We take a hint from \cite{20,20b} and select a PD controller having bidirectional coupling, with gain on position feedback symmetric towards the preceding and following vehicle, but gain on velocity asymmetric. Considering a simple sample-and-hold digital actuation, we will thus assume that $\;u_k(\tau) = $
\begin{eqnarray*}
f(t,k) & := & b_1 (v_{k-1}(t) - v_k(t)) + b_2 (v_{k+1}(t) - v_k(t))\\
&& +  a (x_{k-1}(t) -x_{k}(t)) + a (x_{k+1}(t) - x_k (t))
\end{eqnarray*}
for all $\tau \in (t,t+\dt]$ in continuous-time. After exact integration this yields the exact discrete-time model \eqref{eq:model} with
\begin{eqnarray*}
u_{k,1}(t) & = & f(t,k) \, \dt \; , \quad u_{k,2}(t) \,=\, f(t,k) \, \dt^2/2 \; .
\end{eqnarray*}
The simulation takes arbitrary values $a=1$, $b_1=2$, $b_2=0.5$, and $\dt=0.1$; for the first and last vehicle, we just drop from the feedback law the term associated to the missing neighbor. 

The precise scaling of input disturbances to apply and of output signals to monitor, depends on the definition of string stability that one wishes to consider. We will illustrate the BIBO type scaling of Definition 4, with $\alpha=1$ independent of $N$. As the illustrated controller is linear, it is just a matter of scaling to translate the simulation to other definitions.

\begin{figure}
\includegraphics[width=45mm]{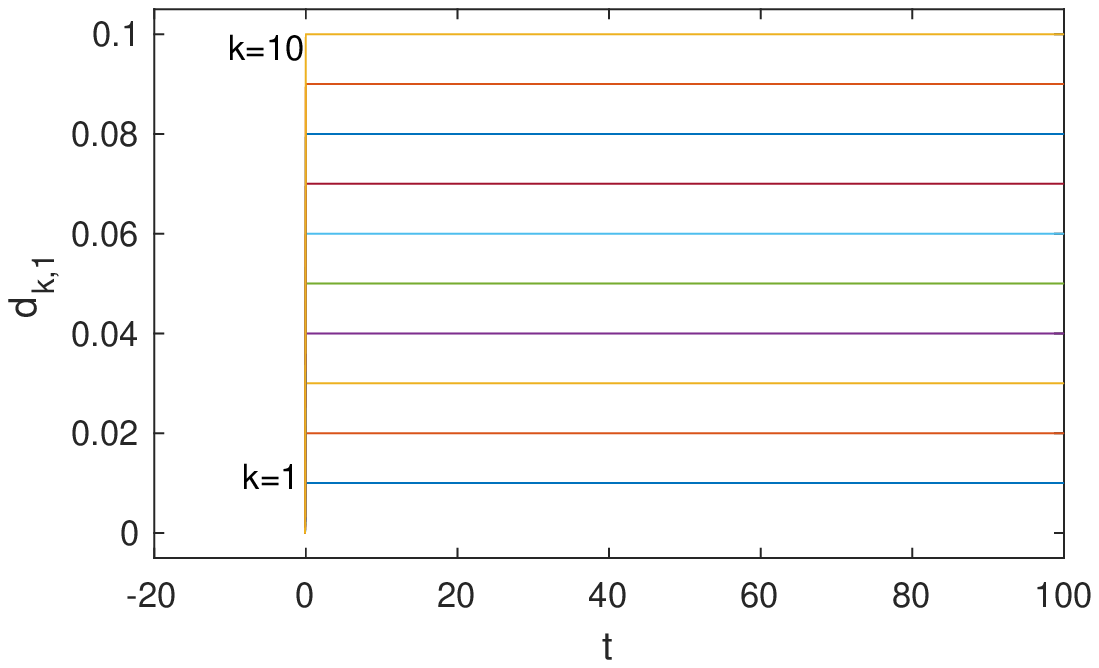}\hspace{-3mm} \includegraphics[width=45mm]{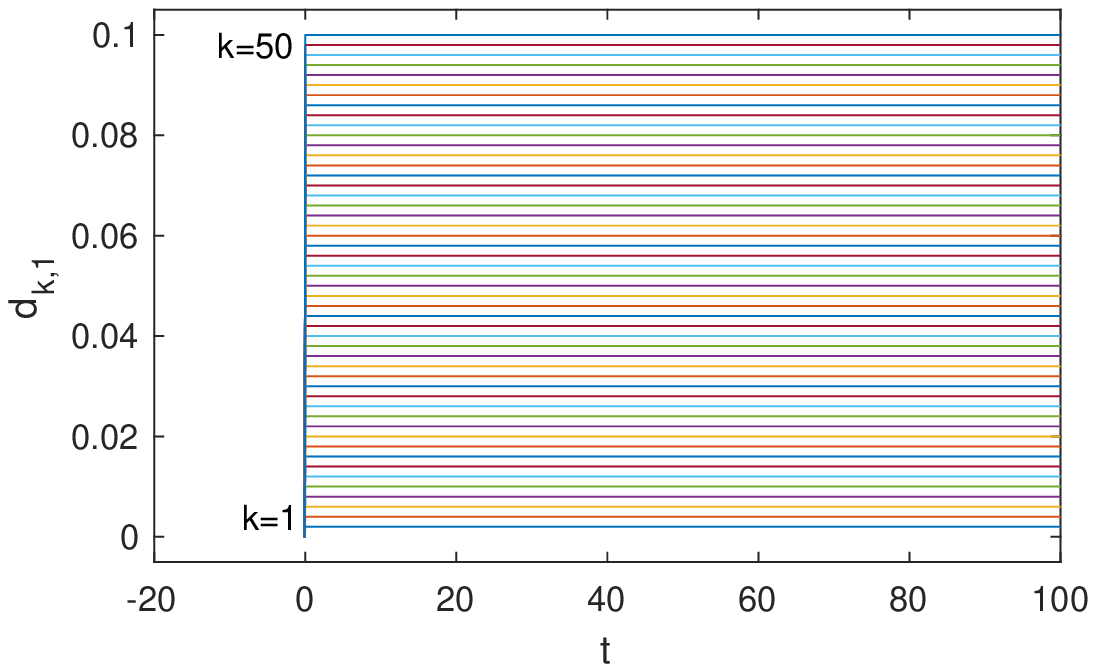}
\caption{Disturbance inputs applied to the vehicles, according to \eqref{eq:disturbance} but without limiting the time window to $t<T$ (see main text), for $N=10$ and for $N=50$ vehicles respectively. The figure is showing $d_{k,1}(t)$ as the corresponding $d_{k,2}(t)$ are just the same multiplied by $\dt$.}\label{fig:input}
\end{figure}

To better illustrate the effect of $N$ only, we will make a small variation on the applied disturbance. Indeed, since our analysis in Lemma 1 only goes up to time $t=T$, for a causal system, it does not matter which disturbance we apply for $t>T$. Therefore, we take an input disturbance of the form \eqref{eq:disturbance}, but instead of stopping it at $t=T$ which scales with $N$, we apply this same nonzero input for all $t>0$, i.e.~we are examining a sort of step response. In this way, the input disturbance acting on vehicle say $k=N/2$ becomes independent of $N$, as illustrated on Fig.~\ref{fig:input}. The impossibility to achieve BIBO string stability via Lemma 1 will be visible as the fact that, when $N$ grows, the maximum deviation $|e_{N/2}(T)|$ at time $T=N\,\dt/5$ grows unboundedly.

\begin{figure}
\includegraphics[width=80mm]{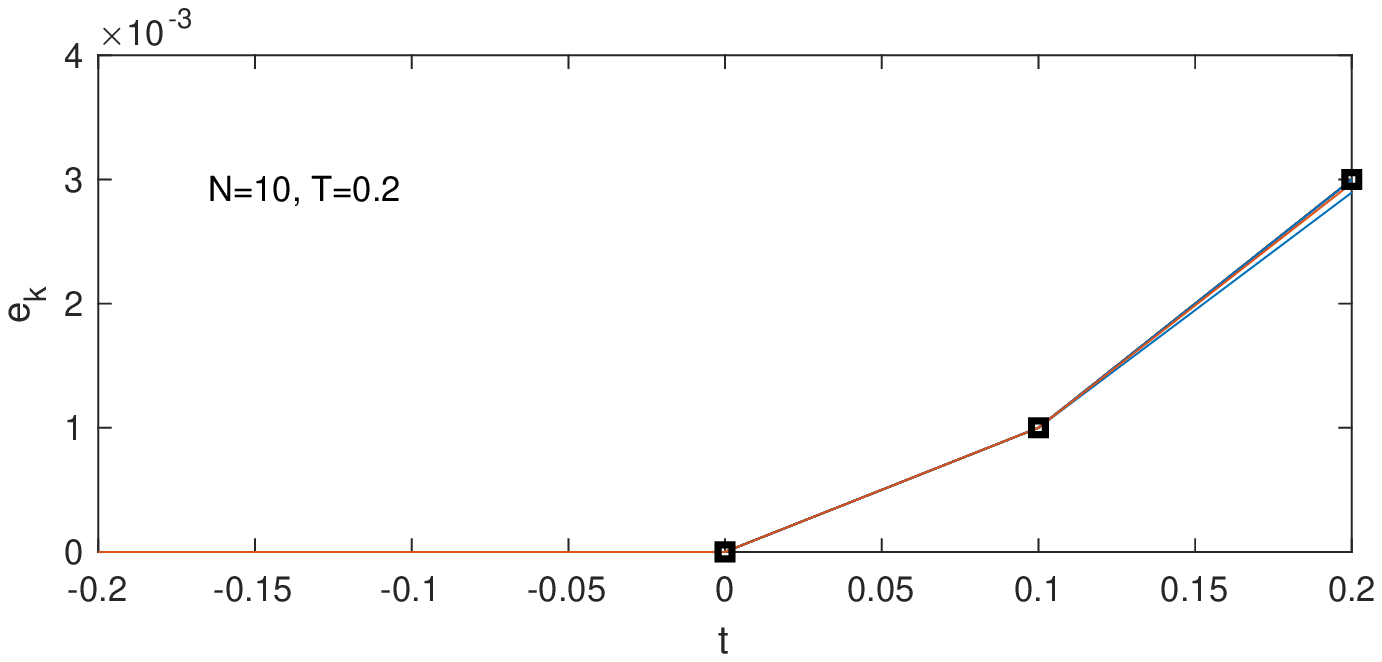} \includegraphics[width=80mm]{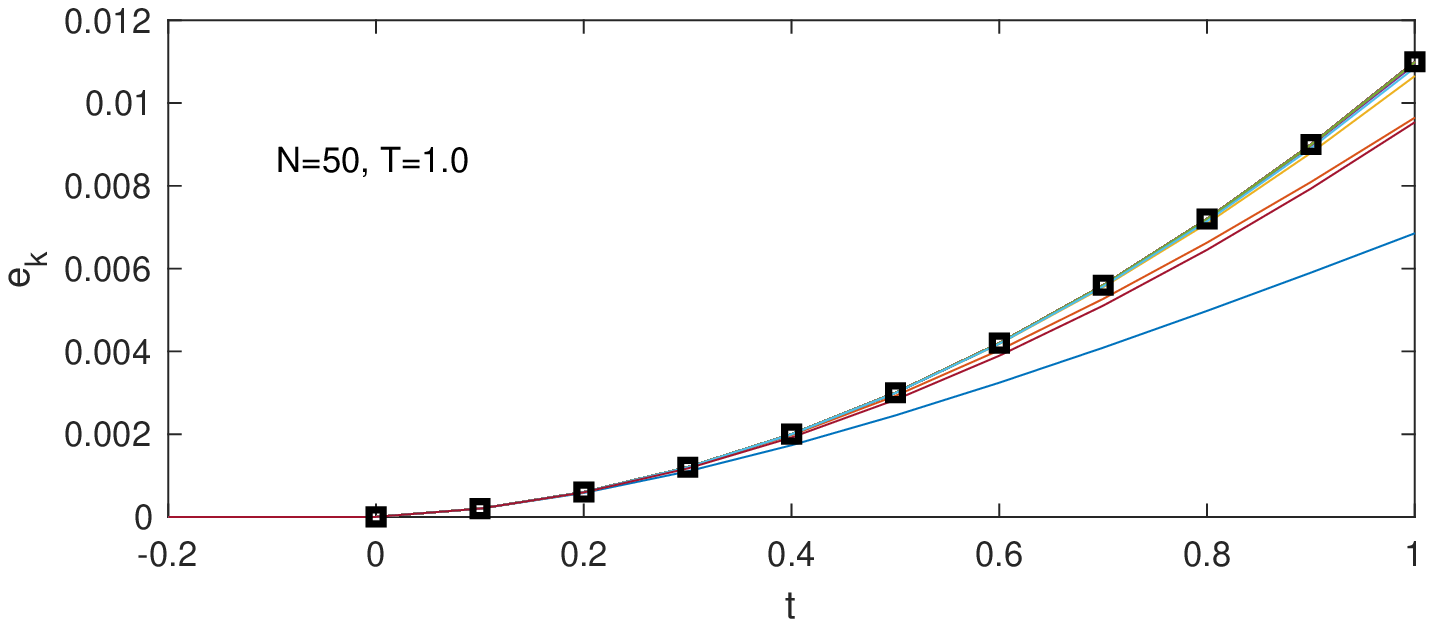}
\caption{Evolution of the distance errors $e_k(t)$ for $t\geq T$, for a vehicle chain (see main text for details) subject to the input disturbances shown on Fig.~\ref{fig:input}, and for $N=10$ and $N=50$ respectively. Note the different scales on \emph{both} axes. In agreement with our analysis we have taken $T=N \dt / 5$. For this time interval, a significant number of $e_k(t)$ is supposed to follow the solution described by Lemma 1; the latter is plotted as black dots which indeed superimpose with a number of simulated curves.}\label{fig:result}
\end{figure}

On Figure \ref{fig:result}, we show the simulated spacing error between consecutive vehicles under this model, up to $t=T$ and for two different chain lengths $N=10$ and $N=50$. The black squares are the solution given by Lemma 1. One indeed observes that many (central) vehicles follow this solution, while others are progressively affected by the boundary effects and behave differently. Note that as $N$ increases, in accordance with Lemma 1, the error at a given time becomes lower; this is due to the consecutive vehicles' disturbance inputs becoming more similar as $N$ increases. However, in return, the solution of Lemma 1 remains valid for a longer time $T$ and therefore overall the error at time $T$, that we can easily compute, keeps increasing unboundedly with $N$.

\begin{figure}
\includegraphics[width=90mm]{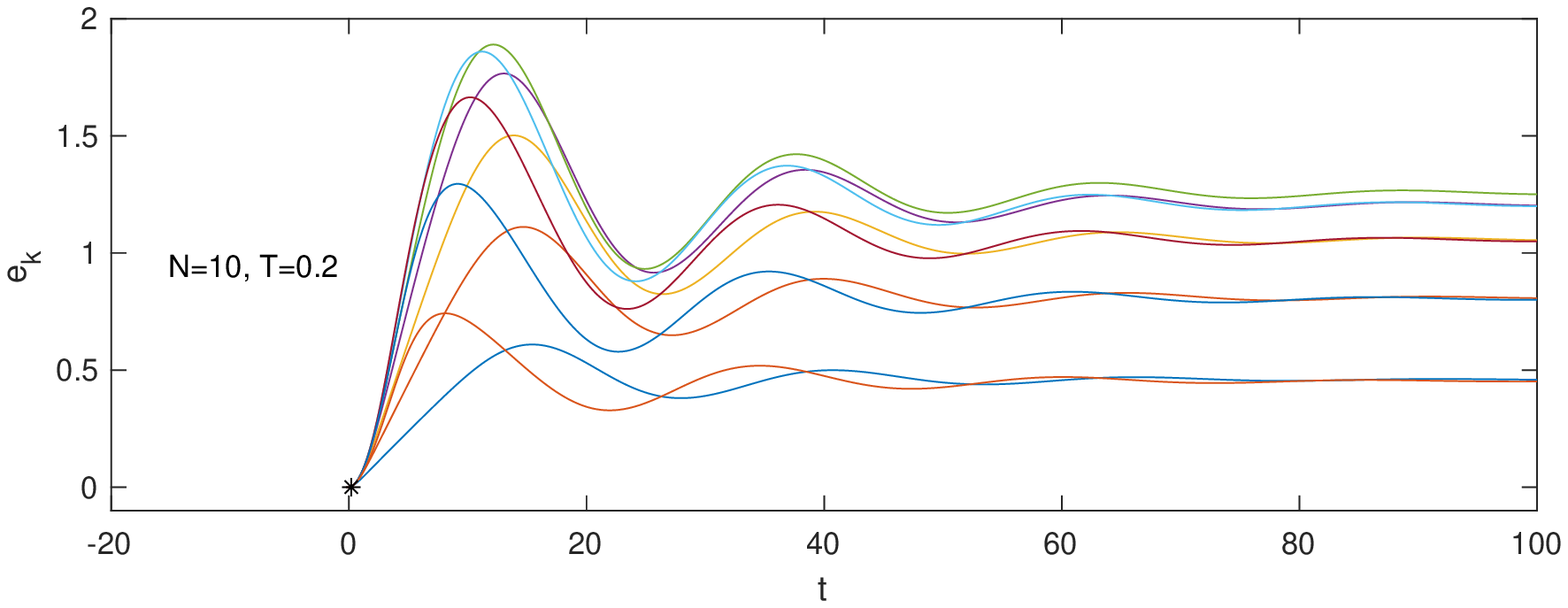} \includegraphics[width=90mm]{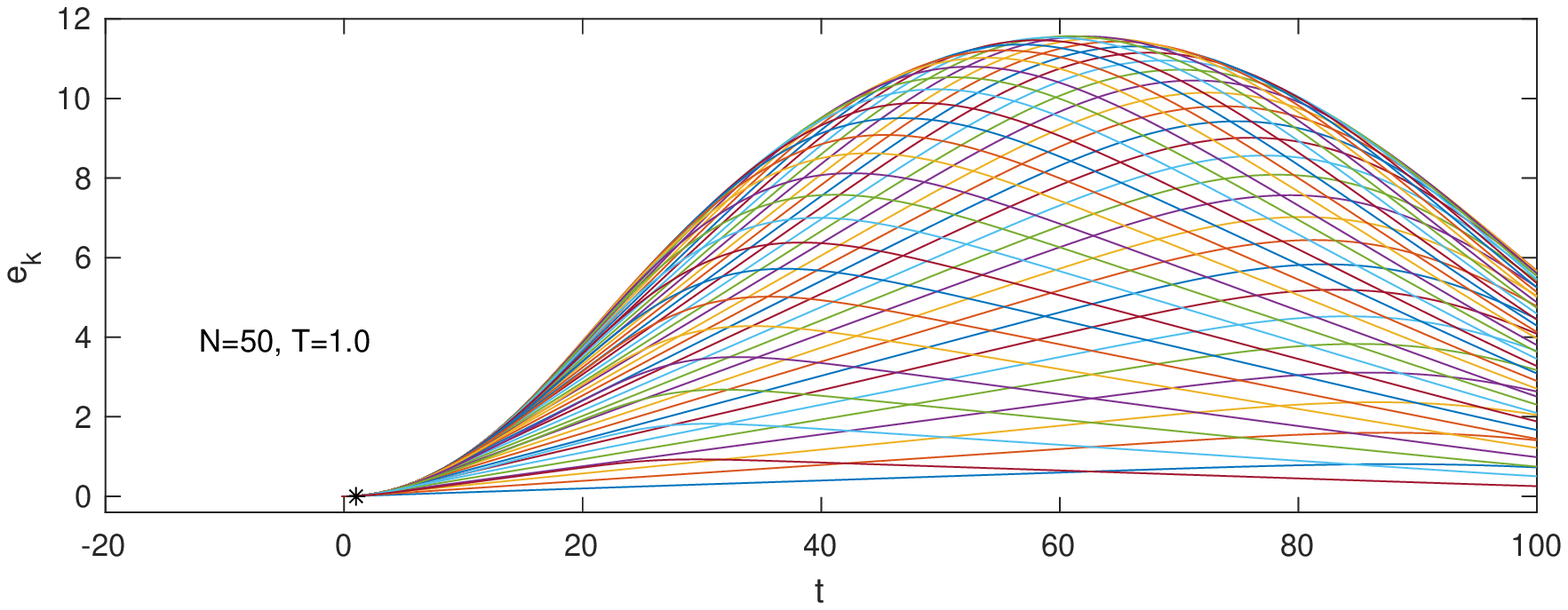}
\caption{Evolution of the $e_k(t)$ for larger times $t$, i.e.~the plots of Figure \ref{fig:result} are in fact zooms on the beginning of the present plots. For $O(t)>T$ (marked with a black dot), the result of Lemma 1 no longer holds and the system behavior does depend on the particular controller choice. We observe that, with the chosen controller: (i) for fixed $N$ indeed each error eventually stabilizes to a bounded value; and (ii) the analysis of Lemma 1, i.e.~with errors taken into account only up to time $t=T$, is very optimistic.}\label{fig:allt}
\end{figure}

We can also have a look at the behavior of the chain for $t>T$, see Figure \ref{fig:allt}. The most striking observation is that the errors keep increasing way beyond the point covered by our analysis (black dot very close to the origin). However, as this depends on the chosen controller we must be careful about further negative conclusions. The main conclusion might thus just be that for a fixed $N$, the chosen controller indeed stabilizes the errors to bounded values, i.e.~it does effectively stabilize the system. Only, it does not so uniformly in $N$, and this is what string instability essentially means.

\section{Concluding discussion}\label{sec:conclusion}

This paper significantly extends the scope of an impossibility result regarding string stability towards input disturbances acting on all subsystems. Indeed, while existing results have centered on LTI systems, we here allow controllers to be nonlinear, $N$-dependent, time-varying, thus possibly modulated and digitally quantized --- as well as using any type of perfect local communication at finite speed. The analysis involves no complicated elements once the setting and example are identified, but as the search for alternative controllers had remained open so far, it appears to give a definite answer clearly narrowing down the options towards achieving string stability. Essential features for the impossibility are:
\begin{itemize}
\item Second-order integrator model for individual subsystems: if the dynamics was first-order, our counterexample would not work;
\item Relative measurements: variations that do solve string stability by adding an absolute velocity term are known, see e.g.~time-headway policies~\cite{10,12,13}. With respect to this criterion, academically, string instability appears more than ever as a property of distributed sensing. In practice, using absolute velocity in the feedback controller or damping becomes a question of hardware and application tradeoff. We must mention that the string instability issue is \emph{not} directly linked to the low observability for long-range modes in distributed systems with relative measurements \cite{MyEELT}. Indeed, here the target variables are not the absolute displacements $x_k$, for which indeed there would be an observability issue, but rather the relative displacements $e_k$, which are directly measured. Also see the previous point.
\item Homogeneous controller, i.e.~same logic with same parameter values at all vehicles: technically, the possibility remains that heterogeneous controllers, i.e.~letting the different vehicles react \emph{differently} to the same signals, could solve the issue. However, we currently have no clue how to design this heterogeneity --- unless one would allow parameters increasing unboundedly with chain length $N$, which however would pose other obvious problems. Controllers periodic in vehicle number, do not seem to work.
\item Discrete-time controller: this should be representative in practice of a realistic digital controller. Rigorously, our counterexample analysis would break down when reducing the discretization step $\dt$ with $N$. However, a property that only holds with infinitely large bandwidth $1/\dt$ for communication and/or control, is usually not robust in practice; this suggests that any ``reasonable'' continuous-time controller would fail too. Note that the standard string stability model here includes no measurement nor communication imperfections, while with extreme continuous-time controllers that are badly modeled by finite $\dt$ those can be expected to become important. 
\end{itemize}
Also note that we have only identified one particular, badly rejected disturbance input. In practice, for a generic disturbance, the situation might often be better, but also worse.\\

With this we believe to have given at least a much more comprehensive picture of what can be done on the standard academic property of string stability. If this string stability property appears critical in some key applications, those results should help guide a possible search for very particular controllers to achieve it, if it is feasible at all without relying on absolute velocities. A point that we did not study is string stability with respect to disturbances on the initial state, instead of on input signals; a similar analysis might be possible.

A different option for the future is to acknowledge that the academic definition of string stability is too strong to be useful, even in an extended framework with nonlinear controllers and so on. In that sense, we can think of two reasonable variations on string stability.
\begin{itemize}
\item One option is to consider the tradeoff in a more integrated picture for finite $N$: to have a given acceptable error, what are the best possible combinations of chain length $N$, absolute-velocity-dependence $h$, control+communication bandwidths $1/\dt$, possibly nonlinear effects, and associated gains in presence of other noises? This, knowing that the limit for infinite $N$ will not work, but will also not be essential for most applications.
\item Another approach would be to acknowledge that the worst-case formulation of string stability is too strong: as the worst-case disturbance could become more and more unlikely with increasing $N$, it might be more telling to take the limit $N\rightarrow \infty$ with a probability distribution over disturbances. In \cite{18}, precisely this approach is taken for the behavior of a lattice of simple \emph{linearly coupled} systems.\\ 
\end{itemize}

As a final word, we may reflect on the more profound implications of our impossibility result. The investigation of \cite{18} for instance is motivated by the stability of physical matter, which after all appears to be governed by forces depending on \emph{relative} states. Implications are also expected for the numerical simulation of related PDEs. It may be an imortant theoretical aim to pin down what essential element in the system structure leads to this impossibility.

A first point in this direction is that the analysis towards our Theorem 2 can be easily extended to other spatial interconnections structure, e.g.~a $D$-dimensional \emph{lattice} of $N$ possibly \emph{nonlinear} systems:
\begin{itemize}
\item We can keep our counterexample with $d_k$ increasing along one dimension of the lattice from $0$ to $\alpha$ with steps $\tfrac{\alpha}{N^{1/D}}$, and constant along the other dimensions. 
\item Computing the acceptable $T$ and $\alpha$ for each case, we get the relevant error growing like $N^{\beta}$ where $\beta\geq0$ depends on the lattice dimension and on the choice of definition, but it is always $>0$ for Definitions 2-4.
\end{itemize}
Compared to \cite{18}, we thus generalize the setting by allowing any nonlinear, time-varying local interactions towards improving the situation, but we obtain a more negative result by considering \emph{the worst} disturbance distribution, over time and over subsystem indices. If your aim is to break a system into parts, this particular disturbance may be useful insight. In contrast, to understand the stability of lattices in a natural environment, one may \emph{have to} acknowledge that bad disturbances in fact become negligibly probable with increasing $N$. This sets a maybe unexpected link between distributed systems and error correcting codes, where scaling to larger codes must essentially rely on the increasing unlikelihood of uncorrectable errors \cite{ErrCorrCodes}.


%
%
%
\end{document}